\documentclass{amsart}
\usepackage{amsmath}
\usepackage{amssymb}
\usepackage{amsthm}

\DeclareMathOperator{\tr}{tr}

\DeclareMathOperator{\dvol}{dvol}

\DeclareMathOperator{\Ric}{Ric}

\DeclareMathOperator{\CD}{CD}

\newcommand{\lp}{\langle}
\newcommand{\rp}{\rangle}
\newcommand{\lv}{\lvert}
\newcommand{\rv}{\rvert}
\newcommand{\lV}{\lVert}
\newcommand{\rV}{\rVert}

\newcommand{\mE}{\mathcal{E}}
\newcommand{\mF}{\mathcal{F}}

\newcommand{\mQ}{\mathcal{Q}}

\newcommand{\mT}{\mathcal{T}}

\newcommand{\mW}{\mathcal{W}}

\newcommand{\kX}{\mathfrak{X}}

\newcommand{\bR}{\mathbb{R}}

\newtheorem{thm}{Theorem}[section]
\newtheorem{prop}[thm]{Proposition}
\newtheorem{lem}[thm]{Lemma}
\newtheorem{cor}[thm]{Corollary}

\theoremstyle{definition}
\newtheorem{defn}[thm]{Definition}

\theoremstyle{remark}
\newtheorem{remark}[thm]{Remark}

\numberwithin{equation}{section}

\begin{document}

\title[Conformal invariants and GNS inequalities]{Conformal invariants measuring the best constants for Gagliardo-Nirenberg-Sobolev inequalities}
\author{Jeffrey S. Case}
\thanks{Partially supported by NSF-DMS Grant No.\ 1004394}
\address{Department of Mathematics \\ Princeton University \\ Princeton, NJ 08544}
\email{jscase@math.princeton.edu}
\keywords{Gagliardo-Nirenberg-Sobolev inequality; sharp constants; smooth metric measure space; conformal geometry}
\subjclass[2000]{Primary 49Q20; Secondary 53C21}
\begin{abstract}
We introduce a family of conformal invariants associated to a smooth metric measure space which generalize the relationship between the Yamabe constant and the best constant for the Sobolev inequality to the best constants for Gagliardo-Nirenberg-Sobolev inequalities $\lV w\rV_q\leq C\lV\nabla w\rV_2^\theta \lV w\rV_p^{1-\theta}$.  These invariants are constructed via a minimization procedure for the weighted scalar curvature functional in the conformal class of a smooth metric measure space.  We then describe critical points which are also critical points for variations in the metric or the measure.  When the measure is assumed to take a special form --- for example, as the volume element of an Einstein metric --- we use this description to show that minimizers of our invariants are only critical for certain values of $p$ and $q$.  In particular, on Euclidean space our result states that either $p=2(q-1)$ or $q=2(p-1)$, giving a new characterization of the GNS inequalities whose sharp constants were computed by Del Pino and Dolbeault.
\end{abstract}
\maketitle


\section{Introduction}
\label{sec:intro}

Understanding the best constants in analytic inequalities plays an important role in geometric analysis and the many fields it interacts with.  Perhaps the most famous example comes via the Sobolev inequality and its importance in the resolution of the Yamabe problem.  The classical Sobolev inequality on the Euclidean space $\bR^n$ states that there is a constant $C_S$ such that
\begin{equation}
\label{eqn:sobolev}
\lV w\rV_{\frac{2n}{n-2}}^2 \leq C_S \lV\nabla w\rV_2^2
\end{equation}
for all $w$ in the Sobolev space $L_1^2(\bR^n)$.  On the other hand, given a compact Riemannian manifold $(M^n,g)$ with $n\geq 3$, the Yamabe problem asks when one can find a metric of constant scalar curvature in the conformal class $[g]$, and is equivalent to finding a smooth function which realizes the Yamabe constant
\begin{equation}
\label{eqn:yamabe}
\sigma_1(g) = \left\{ \frac{\int (\lv\nabla w\rv^2+\frac{n-2}{4(n-1)}Rw^2)\dvol}{(\int w^{\frac{2n}{n-2}}\dvol)^{\frac{n-2}{n}}} \colon 0\not=w\in L_1^2(M) \right\} ,
\end{equation}
where $R$ denotes the scalar curvature of $g$.  By the definition of the Yamabe constant, it is unsurprising that the Sobolev inequality should play an important role in the analysis of the Yamabe problem.  What is surprising is the importance of knowing the sharp constant $C_S$ in~\eqref{eqn:sobolev}: The standard resolution of the Yamabe problem~\cite{Aubin1976,Schoen1984,Trudinger1968,Yamabe1960} requires showing that for a compact manifold not conformally equivalent to the standard sphere, the Yamabe constant is strictly less than $C_S^{-1}$, the Yamabe constant of the standard sphere.

That the sharp constant $C_S$ in the Sobolev inequality is the reciprocal of the Yamabe constant~\eqref{eqn:yamabe} on the sphere is intuitively clear due to its conformal invariance and geometric interpretation; one expects the sharp constant to be realized by the function
\begin{equation}
\label{eqn:sobolev_sharp}
w(x) = \left(1+\lv x\rv^2\right)^{-\frac{n-2}{2}},
\end{equation}
which is the conformal factor relating the Euclidean metric $g_0$ to a constant curvature metric $g=w^{\frac{4}{n-2}}g_0$ on the $n$-sphere.  That this is indeed the case was first proven by Aubin~\cite{Aubin1976s} and Talenti~\cite{Talenti1976}, and since then many different proofs, focusing on different interpretations of the Sobolev constant, have appeared.  Particularly noteworthy from the perspective of this article is the recent rearrangement-free proof given by Frank and Lieb~\cite{FrankLieb2012b}, which relies only on the conformal invariance and as such, leads to a proof of the equivalent inequality in the Heisenberg group~\cite{FrankLieb2012a}.

The connection between sharp Sobolev inequalities and problems in conformal geometry also appears in the consideration of higher order conformally covariant operators.  More precisely, there is a close connection between Sobolev inequalities of the form $\lV w\rV_{\frac{2n}{n-2k}}^2\leq C\lV\nabla^kw\rV_2^2$ and the study of conformally covariant operators with leading order term $(-\Delta)^k$.  Again, specific powers of the functions~\eqref{eqn:sobolev_sharp} can be shown to be extremal functions for these Sobolev inequalities, which is easily seen from the geometric perspective.  For details and references, we refer to the recent survey by Chang~\cite{Chang2005}.

One goal of this article is to pursue the connection between conformal geometry and another family of Sobolev-type inequalities, namely the Gagliardo-Nirenberg-Sobolev (GNS) inequalities
\begin{equation}
\label{eqn:gns}
\lV w\rV_q \leq C \lV\nabla w\rV_2^\theta \, \lV w\rV_p^{1-\theta}
\end{equation}
for $1\leq p\leq q\leq\frac{2n}{n-2}$, where $\theta$ is determined by scaling; i.e.\ \eqref{eqn:gns} is invariant under dilations of $\bR^n$.  This inequality includes as special cases the Sobolev inequality when $q=\frac{2n}{n-2}$ and Gross' logarithmic Sobolev inequality~\cite{Gross1975} in the limiting case $p,q\to2$ (cf.\ \cite{DelPinoDolbeault2002}).

While the sharp constants for the GNS inequalities~\eqref{eqn:gns} are not known in general, they have been computed for two special families of exponents by Del Pino and Dolbeault~\cite{DelPinoDolbeault2002}.  Explicitly, they computed the best constants for the family
\begin{equation}
\label{eqn:gns_dd}
\lV w\rV_{2t} \leq C \lV\nabla w\rV_2^\theta \lV w\rV_{t+1}^{1-\theta}
\end{equation}
with $t\in(1,\frac{n}{n-2}]$, and also for the family
\begin{equation}
\label{eqn:gns_dd_other}
\lV w\rV_{t+1} \leq C \lV\nabla w\rV_2^\theta \lV w\rV_{2t}^{1-\theta}
\end{equation}
with $t\in[\frac{1}{2},1)$; indeed, this also includes in the limiting case $t\to 1$ the best constant in the logarithmic Sobolev inequality.  The extremal functions are
\begin{equation}
\label{eqn:gns_dd_minimizer}
w(x) = (1+\lv x\rv^2)^{-\frac{1}{t-1}}
\end{equation}
for the inequalities~\eqref{eqn:gns_dd} and
\begin{equation}
\label{eqn:gns_dd_minimizer_other}
w(x) = (1-\lv x\rv^2)_+^{-\frac{1}{t-1}}
\end{equation}
for the inequalities~\eqref{eqn:gns_dd_other}.  These functions admit natural interpretations in terms of the conformal equivalence of Euclidean space with the sphere and with hyperbolic space, respectively.  In particular, one is naturally led to wonder if there are Yamabe-type constants on a conformal manifold which are equivalent on Euclidean space to the best constants in the GNS inequalities~\eqref{eqn:gns}, and if these constants can be used to ``explain'' the special nature of the families~\eqref{eqn:gns_dd} and~\eqref{eqn:gns_dd_other}.  We note that a different explanation has recently been given by Agueh~\cite{Agueh2008} using ideas from optimal transportation, based upon the transport-based derivation of the best constants in these inequalities by Cordero-Erausquin, Nazaret and Villani~\cite{CorderoErausquinNazaretVillani2004}.

Due to the three different norms appearing in a general GNS inequality, an important task in formulating the desired Yamabe-type constants is to find a suitable interpretation of these norms.  One way to approach this problem is by using smooth metric measure spaces (cf.\ \cite{CarlenFigalli2011,CorderoErausquinNazaretVillani2004}).  In this way, we will arrive at conformal invariants associated to smooth metric measure spaces which are the desired Yamabe-type constants.  At present, there does not seem to be a consensus on what it means to study a smooth metric measure space as an object in conformal geometry (cf.\ \cite{Case2010a,CGY0}).  For this reason, the second goal of this article is to articulate a geometrically natural formulation of conformal transformations of smooth metric measure spaces, with our conformal invariants as a model for the utility of this perspective.

As a geometric object, a smooth metric measure space is a Riemannian manifold $(M^n,g)$ together with a smooth measure $e^{-\phi}\dvol_g$ --- that is, $\phi\in C^\infty(M)$ and $\dvol_g$ is the Riemannian volume element --- and a dimensional parameter $m\in[0,\infty]$.  The role of the dimensional parameter $m$ is to specify that $e^{-\phi}\dvol_g$ should be regarded as an $(m+n)$-dimensional measure.  As will be made precise in Section~\ref{sec:bg/smms}, this allows one to introduce a curvature tensor $\Ric_\phi^m$ which, via the inequality $\Ric_\phi^m\geq Kg$, characterizes the curvature-dimension bound $\CD(K,m+n)$~\cite{BakryEmery1985}.  Viewed this way, there are two known Yamabe-type constants already defined on smooth metric measure spaces, namely Perelman's $\nu$-entropy~\cite{Perelman1} and the author's $m$-energy~\cite{Case2010b}.

Perelman's $\nu$-entropy can be regarded as a geometric invariant associated to smooth metric measure spaces with $m=\infty$.  From an analytic viewpoint, the $\nu$-entropy is closely related to the logarithmic Sobolev inequality.  A key observation of Perelman is that the $\nu$-entropy is monotone along the Ricci flow, and plays an important role in the study of Type I singularities.  From a geometric viewpoint, this allows one to establish the crucial noncollapsing property of such singularities, in many ways exemplifying the close connections between Sobolev inequalities, logarithmic Sobolev inequalities, and isoperimetric inequalities.

The author's $m$-energy can be regarded as the analogous geometric invariant for arbitrary $m$.  In the simplest case where $\phi=0$ and $m<\infty$, the $m$-energy is the constant
\begin{equation}
\label{eqn:menergy}
\sigma_{1,2} = \inf\left\{ \frac{\left( \int (\lv\nabla w\rv^2 + \frac{m+n-2}{4(m+n-1)}Rw^2)\dvol\right) \left( \int w^2 \right)^{m/n}}{\left( \int w^{\frac{2(m+n)}{m+n-2}}\right)^{(m+n-2)/n}} \right\} ,
\end{equation}
where the infimum is taken over all $0\not=w\in L_1^2(M)$ and $R$ is the scalar curvature of $g$.  These constants are closely related to a different family of GNS inequalities than those studied by Del Pino and Dolbeault~\cite{DelPinoDolbeault2002}, namely
\begin{equation}
\label{eqn:gns_menergy}
\lV w\rV_{\frac{2(m+n)}{m+n-2}} \leq C \lV\nabla w\rV_2^\theta \lV w\rV_2^{1-\theta} .
\end{equation}
Additionally, the constants $\sigma_{1,2}$ are closely related to the Sobolev inequalities one can deduce for smooth metric measure spaces satisfying $\Ric_\phi^m\geq g$; see~\cite{Case2010b} for details.  Moreover, they clearly generalize the Yamabe constant~\eqref{eqn:yamabe}, and can be easily shown to yield Perelman's $\nu$-entropy by differentiating at the endpoint $m=\infty$ (cf.\ \cite{Case2010b,DelPinoDolbeault2002}).  Like the $\nu$-entropy, the constants $\sigma_{1,2}$ were introduced as natural geometric invariants which control local noncollapsing, as needed to address the convergence question considered in~\cite{Case2010b}.

One interesting fact, which is closely related to the difficulty of computing the best constants for general GNS inequalities~\eqref{eqn:gns}, is that the constant $\sigma_{1,2}$ is better behaved on hyperbolic space than Euclidean space.  In the language of~\cite{Case2010b}, this is because hyperbolic space is the only spaceform admitting a nontrivial quasi-Einstein metric.  A slightly different perspective on this will be given in Section~\ref{sec:qe}, where the behavior of GNS inequalities on spaceforms is classified by the sign of the curvature.  A more detailed description of the special nature of the constants $\sigma_{1,2}$ on hyperbolic space, as well as their interpretation in Euclidean space, is contained in forthcoming work of Chang, Yang and the author~\cite{CCY2012}.

In this article, we introduce a family of conformal invariants generalizing the constants $\sigma_{1,2}$ in such a way that they include all of the GNS inequalities~\eqref{eqn:gns} when evaluated on the standard Euclidean space.  Like $\sigma_{1,2}$, our constants are defined as infima of certain natural quotients involving a weighted analogue of the total scalar curvature functional.  We will call these quotients \emph{conformal GNS functionals}, and as will be made precise in Section~\ref{sec:bg/smms} and Section~\ref{sec:functional}, they only depend on the conformal class of a smooth metric measure space.  By studying the variational properties of the conformal GNS functionals, we will then be able to classify those constants for which minimizers are critical points for variations of the metric or the measure, when the measure itself takes a special form; for details, see Section~\ref{sec:qe}.  In particular, restricted to $\bR^n$ with its standard flat metric and Lebesgue measure, our results give the following characterization of the family of GNS inequalities~\eqref{eqn:gns_dd} studied by Del Pino and Dolbeault~\cite{DelPinoDolbeault2002}.

\begin{thm}
\label{thm:main_thm}
The only GNS inequalities~\eqref{eqn:gns} for which extremal functions are also critical points of the conformal GNS functional~\eqref{eqn:gns_functional} for variations in the metric or the measure are the families~\eqref{eqn:gns_dd} and~\eqref{eqn:gns_dd_other}.
\end{thm}

After introducing the aforementioned conformal invariants, the proof of Theorem~\ref{thm:main_thm} hinges upon two observations.  The most important observation is that minimizers which are critical for variations of the metric give, via a pointwise conformal transformation, smooth metric measure spaces for which the Bakry-\'Emery Ricci tensor is a (not necessarily constant) multiple of the metric.  This is a weakening of the usual condition that a smooth metric measure space be quasi-Einstein (cf.\ \cite{Case2011t,Catino2010}), but nevertheless imposes strong restrictions on the underlying smooth metric measure space.

The second observation is that, because of the variational structure of our constants, there are divergence-type formul\ae\ which hold for minimizers of the constants which are critical for variations in the metric or the measure.  In order to give what we feel is a conceptually and notationally simple proof and application of these formul\ae, we will introduce and use some basic aspects of the tractor calculus from conformal geometry~\cite{Bailey1994} to this problem.  For our purposes, the tractor calculus provides a convenient vector bundle in which to perform our computations; for details, see Section~\ref{sec:bg/tractor}.  In particular, we will derive in Section~\ref{sec:lemma} a divergence-type equation which gives a local formulation of Obata's argument classifying constant scalar curvature metrics on compact conformally Einstein manifolds, which we expect to be of independent interest.

Another interesting aspect of the family of GNS inequalities~\eqref{eqn:gns_dd} is that they are stable, in the sense that if a function $\xi$ almost satisfies equality in~\eqref{eqn:gns_dd}, then $\xi$ is close to a function of the form~\eqref{eqn:gns_dd_minimizer}.  For example, Bianchi and Egnell~\cite{BianchiEgnell1991} showed that there is a dimensional constant $\alpha$ such that
\begin{equation}
\label{eqn:bianchi_egnell}
\alpha\inf \lV \nabla(\xi-w)\rV_2 \leq C_S \lV\nabla\xi\rV_2 - \lV\xi\rV_{\frac{2n}{n-2}}
\end{equation}
for all $\xi\in L_1^2(\bR^n)$, where the infimum is taken over all functions $w$ which realize equality in the Sobolev inequality~\eqref{eqn:sobolev}.  Using a nice tensorization argument, Bakry has shown that this stability result can be extended to the family of GNS inequalities~\eqref{eqn:gns_dd} (see~\cite{BakryGentilLedoux2012,CarlenFigalli2011}).  Given that the proof of~\eqref{eqn:bianchi_egnell} depends in an important way on computing the second variation of the Yamabe functional on $\bR^n$, it is natural to wonder if Theorem~\ref{thm:main_thm} states anything about the stability of general GNS inequalities.  We do not treat this question here, but rather leave it as an interesting open question.

This article is organized as follows.

In Section~\ref{sec:bg}, we give the relevant background for our geometric perspective on the GNS inequalities.  First, in Section~\ref{sec:bg/smms}, we explain the principles underlying the study of smooth metric measure spaces as objects in conformal geometry.  Second, in Section~\ref{sec:bg/tractor}, we describe the basic aspects of the tractor calculus as we will use them here.

In Section~\ref{sec:lemma}, we state and prove the aforementioned divergence-type formula localizing Obata's argument, and give a brief discussion on some possible interpretations and generalizations of the result.

In Section~\ref{sec:functional}, we introduce the general Yamabe-type constants on smooth metric measure spaces which, when taken together, characterize the best constants in the GNS inequalities~\eqref{eqn:gns}.  For simplicity, we will initially restrict our attention to GNS inequalities~\eqref{eqn:gns} for which $q\geq 2$.  From their definition, they will easily be seen to share the same basic properties of the Yamabe constant, namely their conformal invariance and their relation to the analogue of the conformal Laplacian on smooth metric measure spaces.

In Section~\ref{sec:qe}, we prove stronger versions of Theorem~\ref{thm:main_thm} for general smooth metric measure spaces.  First, Theorem~\ref{thm:qe_rigid} concerns itself with minimizers of the conformal GNS functionals which are also critical points for variations in the metric.  Second, Theorem~\ref{thm:v_rigid} concerns itself with minimizers which are also critical points for variations in the measure.  Specializing to the case of Euclidean space with its standard measure then recovers Theorem~\ref{thm:main_thm}.

In Section~\ref{sec:conclusion}, we describe the modifications one must make so that our results also include the case $q<2$, and in particular the characterization of the extremal functions~\eqref{eqn:gns_dd_minimizer_other}.

\subsection*{Acknowledgments}
I would like to thank Alice Chang and Paul Yang for our many illuminating conversations on the interplay between conformal geometry and functional inequalities.

\section{Background}
\label{sec:bg}

\subsection{Smooth Metric Measure Spaces and Conformal Geometry}
\label{sec:bg/smms}

In order to make precise the relationship between conformal geometry and the family of GNS inequalities~\eqref{eqn:gns_dd}, we will need to make use of the notion of a smooth metric measure space and discuss in particular conformal transformations of such objects.  To the best of the author's knowledge, there is not a clear consensus as to what the latter idea should mean.  In light of this, we hope that the results of this article can help to clarify this situation (cf.\ \cite{Case2010a,Case2011t,CGY0}).

\begin{defn}
A \emph{smooth metric measure space} is a four-tuple $(M^n,g,e^{-\phi}\dvol,m)$, where $(M^n,g)$ is an oriented Riemannian manifold, $\dvol$ is the corresponding Riemannian volume element, $\phi\in C^\infty(M)$, and $m\in[0,\infty]$ is a dimensional parameter.  In the special case $m=0$, we require $\phi$ to be a constant.
\end{defn}

The role of the dimensional parameter $m$ is to specify that the measure $e^{-\phi}\dvol$ should be regarded as a $(m+n)$-dimensional measure.  This is precisely the sense in which smooth metric measure spaces appear as motivation in the work of Bakry and \'Emery~\cite{BakryEmery1985}, and has led both to a successful generalization of many aspects of comparison geometry to smooth metric measure spaces (e.g.\ \cite{BakryQian2000,Qian1997,Wei_Wylie}) and to a successful synthetic notion of Ricci curvature lower bounds on metric measure spaces~\cite{Lott_Villani,Sturm2006a,Sturm2006b}.

For the purposes of this article, it is enough to know what are the appropriate weighted analogues of the Laplacian, the Ricci curvature, and the scalar curvature on a smooth metric measure space.

\begin{defn}
Let $(M^n,g,e^{-\phi}\dvol,m)$ be a smooth metric measure space.  The \emph{weighted Laplacian} $\Delta_\phi$, the \emph{Bakry-\'Emery Ricci curvature} $\Ric_\phi^m$, and the \emph{weighted scalar curvature} $R_\phi^m$ are given by
\begin{align*}
\Delta_\phi w & = \Delta w - \lp\nabla\phi,\nabla w\rp \\
\Ric_\phi^m & = \Ric + \nabla^2\phi - \frac{1}{m}d\phi\otimes d\phi \\
R_\phi^m & = R + 2\Delta\phi - \frac{m+1}{m}\lv\nabla\phi\rv^2
\end{align*}
for all $w\in C^2(M)$, where our sign convention is $\Delta=\tr\nabla^2$.
\end{defn}

The weighted Laplacian is formally self-adjoint with respect to the given measure, and serves as the natural modification of the usual (rough) Laplacian.  The Bakry-\'Emery Ricci curvature appears when one tries to generalize the Bochner inequality to the weighted setting, and is what gives rise to the characterization of $e^{-\phi}\dvol$ as an $(m+n)$-dimensional measure (cf.\ \cite{BakryEmery1985,Qian1997,Wei_Wylie}).  In this general context, the weighted scalar curvature seems to be relatively new, and seems to have first been introduced in~\cite{Case2010a}.  It is introduced so that the \emph{quasi-Einstein} condition $\Ric_\phi^m=\lambda g$ for a constant $\lambda$ is precisely the Euler-Lagrange equation for the weighted total scalar curvature functional, in analogy with Perelman's introduction of $R_\phi^\infty$ in his celebrated work on the Ricci flow~\cite{Perelman1}.

In order to carry out a meaningful discussion of the conformal geometry of a smooth metric measure space, we must define what it means for two smooth metric measure spaces to be conformally equivalent.

\begin{defn}
Two smooth metric measure spaces $(M^n,g,e^{-\phi}\dvol_{g},m)$ and $(M^n,h,e^{-\psi}\dvol_{h},m)$ are \emph{pointwise conformally equivalent} if there is a function $f\in C^\infty(M)$ such that
\begin{equation}
\label{eqn:conformal_f}
\left( M^n,h,e^{-\psi}\dvol_h,m\right) = \left(M^n,e^{-\frac{2}{m+n-2}f}g,e^{-\frac{m+n}{m+n-2}f}e^{-\phi}\dvol_g,m\right) .
\end{equation}
\end{defn}

In particular, note that the formula for the conformal change of the measure is as the formula for the conformal change of the Riemannian volume element of an $(m+n)$-dimensional manifold, as desired by our interpretation of the dimensional parameter $m$.  Note also that the above definition makes sense in the case $m=\infty$, though now a ``conformal transformation'' is simply a change of measure.

For our purposes, we will be most interested in studying smooth metric measure spaces with $m<\infty$.  In this case, a more intuitive perspective on the rule~\eqref{eqn:conformal_f} is available: If one denotes the measure instead by $v^m\dvol_g$, one has the conformal equivalence
\begin{equation}
\label{eqn:scms_equiv}
\left( M^n,g,v^m\dvol_g, m\right) \sim \left( M^n, e^{2s}g,(e^sv)^m\dvol_{e^{2s}g}, m\right) .
\end{equation}
In other words, $v$ can be naturally regarded as a conformal density (see Section~\ref{sec:bg/tractor}), and we can simply denote a smooth metric measure space by the triple
\begin{equation}
\label{eqn:scms}
(M^n,g,v^m\dvol)
\end{equation}
for $0<v\in C^\infty(M)$, with the dimensional parameter encoded in our notation for the measure.

As we discuss pointwise conformal transformations for smooth metric measure spaces, it will be convenient to introduce some additional notation to succinctly express the effect of a conformal transformation on certain geometric quantities associated to a smooth metric measure space.  In particular, we adapt our notation so that the computations of~\cite{Case2010a,Case2011t} can be used here with minimal reinterpretation.

\begin{defn}
\label{defn:cwm}
Let $(M^n,g,v^m\dvol)$ be a smooth metric measure space and let $u=e^{\frac{f}{m+n-2}}$ be a smooth function on $M$.  Then the \emph{weighted measure $d\omega$}, the \emph{Bakry-\'Emery Ricci curvature $\Ric_{f,\phi}^m$}, and the \emph{weighted scalar curvature $e^{\frac{2f}{m+n-2}}R_{f,\phi}^m$} of the conformally equivalent smooth metric measure space~\eqref{eqn:conformal_f} are given by
\begin{align*}
d\omega & = u^{-m-n}v^m\dvol \\
\Ric_{f,\phi}^m & = R + \nabla^2(f+\phi) + \frac{1}{m+n-2}df\otimes df - \frac{1}{m}d\phi\otimes d\phi + \frac{1}{m+n-2}\Delta_{f+\phi}f\,g \\
R_{f,\phi}^m & = \tr_g\Ric_{f,\phi}^m + m\Delta_{f+\phi}\left(\frac{f}{m+n-2}+\frac{\phi}{m}\right) .
\end{align*}
Additionally, for convenience we introduce the notation
\[ d\rho = u^{2-m-n}v^m\dvol, \qquad \beta = \frac{f}{m+n-2}+\frac{\phi}{m} . \]
In particular, this allows us to write the total weighted scalar curvature of the smooth metric measure space~\eqref{eqn:conformal_f} as
\[ \int_M R_{f,\phi}^m d\rho . \]
\end{defn}

Recall that on a Riemannian manifold, there is an operator $L$, called the conformal Laplacian, which is uniquely defined by requiring that it is conformally covariant, is equal to $-\Delta$ plus a zeroth order term, and $L(1)$ evaluates to a constant multiple of the scalar curvature~\cite{Chang2005}.  Due to these properties, this operator plays an important role in the problem of prescribing scalar curvature within a conformal class; i.e.\ the Yamabe problem.  Unsurprisingly, there is a weighted analogue of the conformal Laplacian which enjoys the same properties.

\begin{defn}
Let $(M^n,g,v^m\dvol)$ be a smooth metric measure space.  The \emph{weighted conformal Laplacian} $L_\phi^m$ is the operator
\[ L_\phi^m = -\Delta_\phi + \frac{m+n-2}{4(m+n-1)}R_\phi^m . \]
\end{defn}

It is straightforward to check (cf.\ \cite{Case2011t}) that this operator is conformally covariant, in the sense that
\begin{equation}
\label{eqn:weighted_conformal_laplacian}
L_\phi^m\left[e^{2s}g,e^{(m+n)s}v^m\dvol\right] = e^{-\frac{m+n+2}{2}s} \circ L_\phi^m\left[g,v^m\dvol\right] \circ e^{\frac{m+n-2}{2}s},
\end{equation}
and thus $L_\phi^m$ has the same qualitative properties as the usual conformal Laplacian.  Indeed, one can use the above formula to immediately deduce that, in the setting of Definition~\ref{defn:cwm} and with $\xi=u^{-\frac{m+n-2}{2}}$,
\[ R_{f,\phi}^m = \frac{4(m+n-1)}{m+n-2}\xi^{-\frac{m+n+2}{m+n-2}} L_\phi^m\xi . \]

\begin{remark}
Rewriting~\eqref{eqn:weighted_conformal_laplacian} with $s=-\frac{f}{m+n-2}$, one realizes that Perelman's operator $L_\phi^\infty=-\Delta_\phi+\frac{1}{4}R_\phi^\infty$ is conformally covariant,
\[ L_\phi^\infty\left[g,e^{-f-\phi}\dvol\right] = e^{\frac{f}{2}}\circ L_\phi^\infty\left[g,e^{-\phi}\dvol\right] \circ e^{-\frac{f}{2}} , \]
providing arguably the most natural perspective on how Perelman's operator fits into the framework of smooth metric measure spaces and conformal geometry (cf.\ \cite{Case2010a,Case2010b,CGY0,CGY2}).
\end{remark}
\subsection{Basic Aspects of the Tractor Calculus}
\label{sec:bg/tractor}

For the remainder of this article, we will only be interested in studying smooth metric measure spaces with $m<\infty$.  In this case, the equivalence relation~\eqref{eqn:scms_equiv} allows one to regard smooth metric measure spaces as objects in conformal geometry.  In particular, by regarding the best constants in the GNS inequalities~\eqref{eqn:gns_dd} as conformal invariants, we can hope to use tools from conformal geometry to better understand them.  For the purposes of this article, the tool we are most interested in is the tractor calculus.  In order to keep our exposition simple, we shall only introduce the tractor calculus as a useful set of vector bundles on a Riemannian manifold, and instead refer the reader to the references~\cite{Bailey1994,BaumJuhl2010,CapSlovak2009} for more formal treatments.

To begin, let $(M^n,c)$ be a conformal manifold; that is, $c=[g]$ for some Riemannian metric $g$, where $h\in[g]$ if and only if there is a smooth function $s$ such that $h=e^{2s}g$.  The \emph{conformal density bundle of weight $w$} is the trivial line bundle whose sections $\sigma\in\mE[w]$ are functions $\sigma\colon c\times M\to\bR$ satisfying
\[ \sigma\left(e^{2s}g,x\right) = e^{ws(x)}\sigma\left(g,x\right) . \]
These spaces give a convenient way to describe functions which change with conformal transformations as powers of the conformal factor, and also provide a convenient way to express conformally covariant operators.  For example, the property~\eqref{eqn:weighted_conformal_laplacian} can be expressed as
\begin{equation}
\label{eqn:weighted_conformal_laplacian_density}
L_\phi^m \colon \mE\left[-\frac{m+n-2}{2}\right] \to \mE\left[-\frac{m+n+2}{2}\right] .
\end{equation}

In one sense, the fundamental object of the tractor calculus is the standard tractor bundle, which is a vector bundle together with a metric and a connection which are canonically associated to a conformal manifold.  For our purposes, we can and will ignore how it is associated to a conformal class.

\begin{defn}
\label{defn:tractor}
Let $(M^n,g)$ be a Riemannian manifold.  The \emph{standard tractor bundle $T$} is the rank $n+2$ vector bundle $\bR\oplus TM\oplus\bR$.  Given a section $I\in\mT=\Gamma(T)$, we will denote
\begin{equation}
\label{eqn:I}
I = \left(\rho \,,\,\omega \,,\, \sigma\right),
\end{equation}
so that $\rho,\sigma\in C^\infty(M)$ and $\omega\in\kX(M)$ is a vector field.

The \emph{tractor metric $h$} is the Lorentzian metric obtained by polarizing
\begin{equation}
\label{eqn:metric}
h(I,I) = 2\sigma\rho + g(\omega,\omega),
\end{equation}
where $I$ is as in~\eqref{eqn:I}.

The \emph{canonical tractor connection $\nabla$} is the connection defined by
\begin{equation}
\label{eqn:connection}
\nabla_x I = \left( \nabla_x\rho-P(\omega) \,,\, \nabla_x\omega+\sigma P(x)+\rho x \,,\, \nabla_x\sigma-g(\omega,x) \right)
\end{equation}
for all $x\in\kX(M)$, where again $I$ is as in~\eqref{eqn:I}, $\nabla$ on the right hand side is the Levi-Civita connection, and $P$ is the Schouten tensor,
\[ P = \frac{1}{n-2}\left(\Ric - \frac{R}{2(n-1)}g\right) . \]

The \emph{splitting operator $L\colon\mE[1]\to\mT$} is the operator
\begin{equation}
\label{eqn:splitting}
L\sigma = \left( -\frac{1}{n}(\Delta\sigma+J\sigma) \,,\, \nabla\sigma \,,\, \sigma \right),
\end{equation}
where $J=\tr P=\frac{R}{2(n-1)}$ is the trace of the Schouten tensor.
\end{defn}

Note in particular that $h$ is a metric connection, $\nabla h=0$.  From a geometric perspective, there are some important implications one can derive from the tractor $Lu$ corresponding to a density $u\in\mE[1]$.

\begin{prop}[cf.\ \cite{Gover2008}]
\label{prop:tractor_properties}
Let $(M^n,g)$ be a Riemannian manifold with standard tractor bundle $T$, and let $u\in\mE[1]$ be positive.  Then
\begin{enumerate}
\item The scalar curvature $R(u^{-2}g)$ of the metric $u^{-2}g$ is given by
\[ R(u^{-2}g) = -n(n-1)\lv Lu\rv^2 , \]
where we have denoted by $\lv\,\cdot\,\rv$ the norm corresponding to the tractor metric $h$, $\lv I\rv^2=h(I,I)$.  In particular, $R(u^{-2}g)$ has constant scalar curvature if and only if $\lv Lu\rv^2$ is constant.
\item The metric $u^{-2}g$ is Einstein if and only if $\nabla Lu=0$.
\end{enumerate}
\end{prop}

\begin{remark}

These results are actually local.  In particular, if one removes the requirement that $u$ be positive, the results still hold for the metric $u^{-2}g$ wherever it is defined; i.e.\ away from the zeroes of $u$.  Thus one can talk about ``almost Einstein metrics'' as densities $u\in\mE[1]$ such that $\nabla Lu=0$, as in~\cite{Gover2008}.  This perspective appears in understanding the minimizers~\eqref{eqn:gns_dd_minimizer_other} of the family~\eqref{eqn:gns_dd_other} of GNS inequalities.
\end{remark}

The second claim of Proposition~\ref{prop:tractor_properties} is the one we are most interested in, and in many ways exemplifies the importance of the tractor calculus to conformal geometry.  First, the $T^\ast M\otimes TM$ component of $\nabla Lu$ according to the definition~\eqref{eqn:connection} is
\begin{equation}
\label{eqn:pde}
\left( uP + \nabla^2 u\right)_0 = \frac{1}{n-2}\left(\Ric(u^{-2}g)\right)_0,
\end{equation}
where $S_0$ denotes the traceless part of a symmetric $(0,2)$-tensor and $\Ric(u^{-2}g)$ is the Ricci curvature of the metric $u^{-2}g$.  This makes it obvious that if $\nabla Lu=0$, then $u^{-2}g$ is Einstein.  The converse of this statement is typically interpreted as the statement that the standard tractor connection prolongs the overdetermined PDE corresponding to finding $u$ such that~\eqref{eqn:pde} vanishes.  Denoting by $X$ the useful tractor $(1,0,0)\in\mT\otimes\mE[1]$, this is equivalent to the following useful lemma.

\begin{lem}
\label{lem:prolongation}
Let $(M^n,g)$ be a Riemannian manifold and let $I\in\mT$ be such that $\nabla I=\alpha\otimes X$ for some one-form $\alpha$.  Then $\alpha=0$.
\end{lem}

Another useful benefit of Proposition~\ref{prop:tractor_properties} is that it implies that the space of (almost) Einstein metrics in a conformal class is a finite dimensional vector space.  Indeed, in Euclidean space, one can give a basis for the standard tractor bundle in terms of orthonormal parallel tractors.  Explicitly, letting $x^1,\dotsc,x^n$ denote the standard Cartesian coordinates on $\bR^n$ and letting $r^2=\lv x\rv^2$ denote the square-distance from the origin, the set
\begin{equation}
\label{eqn:rn_tractor_basis}
\left\{ L\big(\frac{1+r^2}{2}\big), L(x^1), \dotsc, L(x^n), L\big(\frac{1-r^2}{2}\big) \right\}
\end{equation}
of tractors forms a parallel orthonormal basis, with the first tractor having length $-1$ and the rest having length $1$.  Using Proposition~\ref{prop:tractor_properties}, this allows one to easily realize $S^n$ and $H^n$ as conformally equivalent to $\bR^n$, and will allow us to easily deduce the form of the minimizers of the GNS inequalities~\eqref{eqn:gns_dd} and~\eqref{eqn:gns_dd_other} discovered by Del Pino and Dolbeault; see Remark~\ref{rk:dd_proof}.

Finally, to translate between the tensorial expressions for the Bakry-\'Emery Ricci tensor and the weighted scalar curvature of a smooth metric measure space and the tractorial expressions which we will use in our proofs, the following computation from~\cite{Case2011t} will be useful.

\begin{prop}
\label{prop:smms_tractor_formulae}
Let $(M^n,g,v^m\dvol)$ be a smooth metric measure space and let $u=e^{\frac{f}{m+n-2}}$ determine a pointwise conformal change as in~\eqref{eqn:conformal_f}.  Then it holds that
\begin{align*}
\left(\Ric_{f,\phi}^m\right)_0 & = (m+n-2)v\nabla(Lu) - mu\nabla(Lv) \\
& \quad - \frac{1}{n}\left((m+n-2)\lp\nabla(Lu),Lv\rp-m\lp\nabla(Lv),Lu\rp\right)X \\
R_{f,\phi}^m-m\Delta_\rho\beta & = -(m+n-1)nu^{-2}\lv Lu\rv^2 + mn(uv)^{-1}\lp Lu,Lv\rp \\
R_{f,\phi}^m-(m+n)\Delta_\rho\beta & = -(m+n-2)n(uv)^{-1}\lp Lu,Lv\rp + (m-1)nv^{-2}\lv Lv\rv^2 .
\end{align*}
\end{prop}
\section{A Useful Tractor Lemma}
\label{sec:lemma}

A key idea of the proof of Theorem~\ref{thm:main_thm} is to use the corresponding Euler-Lagrange equations to show that a minimizer must be of the form~\eqref{eqn:gns_dd_minimizer}, from which it follows that the GNS inequality in question is of the form~\eqref{eqn:gns_dd}.  It turns out that this observation can be generalized to arbitrary smooth metric measure spaces for which the measure has been specially chosen, as we will show in Theorem~\ref{thm:v_rigid}.  A straightforward proof of this fact can be made by making use of the following lemma.

\begin{lem}
\label{lem:tractor_lemma}
Let $(M^n,g)$ be a Riemannian manifold such that there exists a parallel tractor $Lv\in\mT$.  Then given any tractor $Lu\in\mT$, at each point $p\in M$ for which $u(p),v(p)\not=0$, it holds that
\begin{equation}
\label{eqn:tractor_lemma}
u^{n-2}\delta\left(u^{2-n}\nabla\lp Lu,Lv\rp\right) - \frac{v^{n-1}}{2u}\delta\left(v^{2-n}\nabla\lv Lu\rv^2\right) = -\frac{v}{u}\lv\nabla Lu\rv^2 ,
\end{equation}
where $\delta Y=\tr (Z\mapsto \nabla_ZY)$ is the divergence of a vector field $Y$.
\end{lem}

\begin{proof}

To begin, denote $x=-\frac{1}{n}(\Delta u+Ju)$, so that $Lu=(x,\nabla u,u)$.  By definition, it holds that
\begin{equation}
\label{eqn:dlu2}
\frac{1}{2}\nabla\lv Lu\rv^2 = (uP+\nabla^2u+xg)(\nabla u) + u(\nabla x-P(\nabla u)) = \nabla_{\nabla u}\nabla u + x\nabla u + u\nabla x .
\end{equation}

First, using that $Lv$ is parallel, we compute that
\begin{equation}
\label{eqn:dlulv}
\begin{split}
\nabla\lp Lu,Lv\rp & = \lp\nabla Lu,Lv\rp = (uP+\nabla^2+xg)(\nabla v) + v(\nabla x-P(\nabla u)) \\
& = (uP+\nabla^2+xg)(\nabla v-\frac{v}{u}\nabla u) + \frac{v}{2u}\nabla\lv Lu\rv^2,
\end{split}
\end{equation}
where the last equality follows from~\eqref{eqn:dlu2}.

Next, it is straightforward to check that
\[ \delta(uP+\nabla^2u+xg) = -(n-1)(\nabla x-P(\nabla u)) . \]
Taking the divergence of~\eqref{eqn:dlulv}, it then follows that
\begin{align*}
\Delta\lp Lu,Lv\rp & = \lp uP+\nabla^2 u+xg, \nabla^2v-\frac{v}{u}\nabla^2u\rp + \frac{v}{2u}\Delta\lv Lu\rv^2 + \frac{1}{2u}\lp\nabla\lv Lu\rv^2,\nabla v-\frac{v}{u}\nabla u\rp \\
& \quad - \lp\nabla v-\frac{v}{u}\nabla u, (uP+\nabla^2u+xg)(u^{-1}\nabla u)+(n-1)(\nabla x-P(\nabla u))\rp \\
& = -\frac{v}{u}\lv uP+\nabla^2 u+xg\rv^2 + \frac{n-2}{u}\lp\nabla v-\frac{v}{u}\nabla u,(uP+\nabla^2u+xg)(\nabla u)\rp \\
& \quad + \frac{v}{2u}\Delta\lv Lu\rv^2 - \frac{n-2}{2u}\lp\nabla v-\frac{v}{u}\nabla u,\nabla\lv Lu\rv^2\rp
\end{align*}
where the second equality uses~\eqref{eqn:dlu2} and the fact that $Lv$ is parallel via the general identity
\[ \left(v(uP+\nabla^2 u+xg)-u(vP+\nabla^2v-\frac{1}{n}(\Delta v+Jv)g\right)_0 = \left(v\nabla^2u-u\nabla^2v\right)_0 . \]
Using~\eqref{eqn:dlu2} to eliminate the second summand in the second equality and rewriting the equation to make clear the divergence structure yields~\eqref{eqn:tractor_lemma}.
\end{proof}

It is instructive to reformulate Lemma~\ref{lem:tractor_lemma} in tensorial language.

\begin{cor}
\label{cor:tractor_lemma}
Let $(M^n,g)$ be an Einstein manifold satisfying $\Ric=(n-1)\lambda g$, let $u\in C^\infty(M)$ be a positive function, and let $\hat R$ denote the scalar curvature of the metric $\hat g=u^{-2}g$.  Then, in terms of the metric $g$, it holds that
\[ \frac{1}{(n-2)^2} \left|\Ric(\hat g)_0\right|^2 = \frac{1}{n}u^{n-1}\delta\left(u^{2-n}(\Delta u+n\lambda u)\right) -\frac{n(n-1)}{2}\Delta\hat R . \]
\end{cor}

\begin{proof}

By Proposition~\ref{prop:tractor_properties}, the assumption that $(M^n,g)$ is Einstein implies that $L(1)$ is parallel.  Applying Lemma~\ref{lem:tractor_lemma} together with the formul\ae\ in Definition~\ref{defn:tractor} then yields the result.
\end{proof}

For our purposes, the most important observation about Lemma~\ref{lem:tractor_lemma} is that $\lv\nabla Lu\rv^2$ is nonnegative and vanishes if and only if $u^{-2}g$ is an Einstein metric.  Thus, if one is in a situation where the left hand side of~\eqref{eqn:tractor_lemma} can be written as the divergence of an integrable vector field, one immediately concludes that $Lu$ is parallel.

As is made clear by Corollary~\ref{cor:tractor_lemma}, Lemma~\ref{lem:tractor_lemma} is a slight generalization of the local identity used by Obata to classify on compact manifolds constant scalar curvature metrics which are conformally Einstein~\cite{Obata1971}.  A potentially interesting observation, which we will not prove here, is that Lemma~\ref{lem:tractor_lemma} is easily generalized to the weighted tractor bundles introduced by the author~\cite{Case2011t}.  Moreover, the equivalent tensor formulations makes sense in the limit $m\to\infty$, where it is the elliptic version of Perelman's local monotonicity formula for his $\mW$-functional~\cite{Perelman1}.

Given that Lemma~\ref{lem:tractor_lemma} generalizes to weighted tractor bundles, one might also hope that Lemma~\ref{lem:tractor_lemma} generalizes to other parabolic geometries, or at least other $|1|$-graded parabolic geometries (see~\cite{CapSlovak2009} for examples).  Evidence that this should be the case is found in the generalization by Jerison and Lee~\cite{JerisonLee1988} of Obata's theorem to the setting of CR geometry.  In particular, this might lead to a new derivation of the identity used by Jerison and Lee to make this generalization.
\section{The Conformal Gagliardo-Nirenberg-Sobolev Functional}
\label{sec:functional}

Let us now define the conformal invariants which give the sharp constants in GNS inequalities~\eqref{eqn:gns} when restricted to the standard Euclidean space.

\begin{defn}
Let $(M^n,g,v^m\dvol)$ be a smooth metric measure space and fix $k\in(0,\frac{m+n+2}{2}]$.  The \emph{conformal Gagliardo-Nirenberg-Sobolev functional} $\mQ_k[g,v^m\dvol]$ is given by
\begin{equation}
\label{eqn:gns_functional}
\mQ_k[g,v^m\dvol](w) = \frac{\mF(w) \left(\int_M w^{\frac{2(m+n-k)}{m+n-2}}v^{m-k}\right)^p}{\left(\int w^{\frac{2(m+n)}{m+n-2}}v^m\dvol\right)^q} ,
\end{equation}
where $\mF$ is the energy functional corresponding to the weighted conformal Laplacian,
\[ \mF(w) = \int_M w L_\phi^m w\,v^m\dvol , \]
$p=\frac{2m}{nk}$, and $q=\frac{2m+k(n-2)}{nk}$.
\end{defn}

Note that, if $w=e^{-f/2}$, then
\begin{equation}
\label{eqn:energy_to_curvature}
\mF(w) = \frac{m+n-2}{4(m+n-1)}\int_M R_{f,\phi}^m d\rho .
\end{equation}

By the conformal covariance of the weighted conformal Laplacian~\eqref{eqn:weighted_conformal_laplacian}, it is clear that the conformal GNS functional is conformally invariant, in the sense that
\[ \mQ_k\left[ e^{2s}g, e^{(m+n)s}v^m\dvol_g\right]\left(e^{-\frac{m+n-2}{2}s}w\right) = \mQ_k\left[g,v^m\dvol_g\right]\left(w\right) . \]
This is most easily checked by regarding $w$ as a conformal density, $w\in\mE\left[-\frac{m+n-2}{2}\right]$, and noting that each integrand (including the measure) has total conformal weight zero.

Like the relationship between the Yamabe functional and the sharp Sobolev constant (corresponding above to $m=0$), the relationship between the conformal GNS functionals and the sharp constant in the GNS inequalities~\eqref{eqn:gns} are made by finding the infimum of the conformal GNS functional over all nonvanishing functions $w\in L_1^2(M)$.

\begin{defn}
Let $(M^n,g,v^m\dvol)$ be a smooth metric measure space and fix $k\in(0,\frac{m+n+2}{2}]$.  The \emph{GNS constant} $\sigma_{1,k}=\sigma_{1,k}(g,v^m\dvol)$ is
\[ \sigma_{1,k} = \inf_{0\not=w\in L_1^2(M)} \mQ_k(w) . \]
\end{defn}

Note that, if one takes $M$ to be noncompact and $v$ or $R_\phi^m$ to be unbounded, one must be careful about what is meant by ``$w\in L_1^2(M)$.''  For the purposes of this article, we will only consider the cases where $v$ and $R_\phi^m$ are bounded, so that this becomes a nonissue.  In particular, on the standard Euclidean space, the above definition is an equivalent formulation of the best constant in a GNS inequality.

\begin{lem}
\label{lem:gns_const}
Let $(\bR^n,g,1^m\dvol)$ be the standard Euclidean space for some $m\in[0,\infty)$ and $k\in(0,\frac{m+n+2}{2}]$.  Then, with $\sigma_{1,k}$ the GNS constant of $(\bR^n,g,1^m\dvol)$, it holds that
\[ \lV w\rV_q \leq \sigma_{1,k}^{-\frac{(m+n-2)nk}{(m+n)(2m+k(n-2))}} \lV\nabla w\rV_2^\theta \, \lV w\rV_p^{1-\theta} \]
holds with $q=\frac{2(m+n)}{m+n-2}$ and $p=\frac{2(m+n-k)}{m+n-2}$ and $\theta$ determined by scaling.  In particular, the special family~\eqref{eqn:gns_dd} of GNS inequalities studied by Del Pino and Dolbeault corresponds to the case $k=1$.
\end{lem}

\begin{proof}

This follows immediately from the definition of the conformal GNS functional and the GNS constant $\sigma_{1,k}$.
\end{proof}

In order to prove Theorem~\ref{thm:main_thm}, as well as its more general forms stated in Section~\ref{sec:qe}, the following variational formul\ae\ for the conformal GNS functionals will be useful.

\begin{prop}
\label{prop:variational}
Let $(M^n,g,v^m\dvol)$ be a smooth metric measure space and fix $k\in(0,\frac{m+n+2}{2}]$.  Suppose that $w=e^{-f/2}$ is a critical point of the functional $w\mapsto \mQ_k[g,v^m\dvol](w)$.  Then $w$ satisfies
\begin{equation}
\label{eqn:kcsc}
R_{f,\phi}^m + 2(m+n-k)\mu u^{k-2} v^{-k} = (m+n)\lambda u^{-2},
\end{equation}
where
\begin{equation}
\label{eqn:knorm}
\mu = \frac{m\int R_{f,\phi}^md\rho}{(m+n-2)nk\int (u/v)^kd\omega}, \quad \lambda = \frac{(2m+k(n-2))\int R_{f,\phi}^md\rho}{(m+n-2)nk\omega(M)} .
\end{equation}

If it also holds that $(w,v)$ is a critical point of the functional $(w,v)\mapsto\mQ_k[g,v^m\dvol](w)$, then it holds that
\begin{equation}
\label{eqn:kmu}
R_{f,\phi}^m - m\Delta_\rho\beta + n(2-k)\mu u^{k-2} v^{-k} = n\lambda u^{-2} .
\end{equation}

If instead $(g,w)$ is a critical point of the functional $(g,w)\mapsto\mQ_k[g,v^m\dvol](w)$, then it holds that
\begin{equation}
\label{eqn:kqe}
\Ric_{f,\phi}^m + (2-k)\mu u^{k-2}v^{-k}g = \lambda u^{-2} g .
\end{equation}
\end{prop}

\begin{proof}

For convenience, denote $\mF = \int R_{f,\phi}^m,d\rho$,
\[ \omega_0 = \int_M w^{\frac{2(m+n)}{m+n-2}}v^m\dvol , \qquad \omega_k = \int_M w^{\frac{2(m+n-k)}{m+n-2}}v^{m-k}\dvol, \]
and consider a compactly-supported variation $(\delta g,\delta f,\delta\phi)$ of $(g,f,\phi)$.  Using~\eqref{eqn:energy_to_curvature}, we may thus write
\[ \delta\mQ_k = \frac{(m+n-2)\omega_k^p}{4(m+n-1)\omega_0^q}\left(\delta\mF + \frac{p\mF}{\omega_k}\delta\omega_k - \frac{q\mF}{\omega_0}\delta\omega_0\right) . \]
On the one hand, it is straightforward to check that
\begin{align*}
\delta\omega_0 & = -\int_M \left[ -\frac{1}{2}\tr\delta g + \frac{m+n}{m+n-2}\delta f + \delta\phi \right] w^{\frac{2(m+n)}{m+n-k}} v^m\dvol \\
\delta\omega_k & = -\int_M \left[ -\frac{1}{2}\tr\delta g + \frac{m+n-k}{m+n-2}\delta f + \frac{m-k}{m}\delta\phi \right] w^{\frac{2(m+n-k)}{m+n-2}} v^{m-k}\dvol .
\end{align*}
On the other hand, it is shown in~\cite[Proposition~4.18]{Case2010a} that
\[ \delta\mF = -\int_M \left[ \lp\Ric_{f,\phi}^m-\frac{1}{2}R_{f,\phi}^mg,\delta g\rp + R_{f,\phi}^m \delta f + (R_{f,\phi}^m-2\Delta_\rho\beta)\delta\phi\right] d\rho . \]
Combining these formul\ae\ yields the claimed Euler-Lagrange equations.
\end{proof}

\begin{remark}

As a cautionary note, our convention in defining $\mu$ in~\eqref{eqn:knorm} differs from the more usual convention in the quasi-Einstein literature (cf.\ \cite{Case2010a}).
\end{remark}

\begin{remark}

The cases $m=1$ and $k=1$ of~\eqref{eqn:kmu} and~\eqref{eqn:kqe} give examples of a different type of critical metric which has recently been considered by Miao and Tam~\cite{MiaoTam2008,MiaoTam2011}.  Namely, if~\eqref{eqn:kmu} and~\eqref{eqn:kqe} hold, then it is straightforward to check that the smooth metric measure space $(M^n,u^{-2}g,u^{-n-1}v\,\dvol_g,1)$ satisfies
\[ v\Ric - \nabla^2 v + \Delta v\,g = (n-1)\mu\,g . \]
With $\mu=\frac{1}{1-n}$, this is precisely the equation satisfied by a critical point of the volume functional restricted to the class of constant scalar curvature metrics with prescribed boundary data on a compact manifold with boundary; see~\cite{MiaoTam2008}.
\end{remark}

Let us conclude this section with two useful observations.  First, on compact smooth metric measure spaces, minimizers of the GNS constants always exist.

\begin{prop}
Let $(M^n,g,v^m\dvol)$ be a compact smooth metric measure space with $m\in[0,\infty)$, and fix $k\in(0,\frac{m+n+2}{2}]$.  Then there exists a positive function $w\in C^\infty(M)$ such that $\sigma_{1.k}(g,v^m\dvol)=\mQ_k(w)$.
\end{prop}

\begin{proof}

This is a well-known consequence of the resolution of the Yamabe problem~\cite{Aubin1976,Schoen1984,Trudinger1968,Yamabe1960} in the case $m=0$ and the fact that when $m>0$, the Euler-Lagrange equation~\eqref{eqn:kcsc} is equivalently
\[ \frac{m+n-2}{4(m+n-1)}L_\phi^m w + 2(m+n-k)\mu v^{-k} w^{\frac{m+n+2-2k}{m+n-2}} = (m+n)\lambda w^{\frac{m+n+2}{m+n-2}} ; \]
in particular, it has subcritical Sobolev exponent, and is thus easily solved using a standard minimization argument (cf.\ \cite{DelPinoDolbeault2002}).
\end{proof}

Second, if a minimizer of a GNS constant is also critical for variations in the metric, then it is critical for variations in the measure.

\begin{lem}
\label{lem:critical_g}
Let $(M^n,g,v^m\dvol)$ be a smooth metric measure space and let $w\in C^\infty(M)$ be a minimizer of the GNS constant $\sigma_{1,k}(g,v^m\dvol)$.  If $w$ is also critical for variations in the metric, then it is critical for variations in the measure.
\end{lem}

\begin{proof}

Let $\xi,\psi$ be two smooth, compactly-supported functions, defining the variation $(g,e^{t\xi}u,e^{t\psi}v)$, where $u=w^{-\frac{2}{m+n-2}}$.  By the conformal invariance of the conformal GNS functional and the conformal equivalence
\[ \left(g,e^{t\xi}u,e^{t\psi}v\right) \sim \left(e^{-2t\psi}g,e^{t(\xi-\psi)}u,v\right), \]
we see that $\xi$ and $\psi$ determine a variation of $(g,w)$.  Hence
\[ \frac{d}{dt}\mQ_k\left[g,(e^{-t\xi}v)^m\dvol_g\right]\left((e^{t\psi}u)^{-\frac{m+n-2}{2}}\right) = 0, \]
as desired.
\end{proof}
\section{Solutions Which Are Critical in More Than One Way}
\label{sec:qe}

Let us now turn to proving Theorem~\ref{thm:main_thm} by considering more generally minimizers of the constants $\sigma_{1,k}(g,v^m\dvol)$ which are also critical for variations in $g$ or $v$.  By using the tractor calculus, we will arrive at simple proofs of our results.  To that end, the following reformulation of Proposition~\ref{prop:variational} will be useful.

\begin{prop}
\label{prop:tractor}
Let $(M^n,g,v^m\dvol)$ be a smooth metric measure space, and suppose that $w=u^{-\frac{m+n-2}{2}}$ is a critical point of the conformal GNS functional for variations in $w$.
\begin{enumerate}
\item If also $w$ is a critical point of the conformal GNS functional for variations in $v$, then it holds that
\begin{align}
\label{eqn:lambda} (m+n-2)k\lambda v^2 & = -k(m+n-1)(m+n-2)v^2\lv Lu\rv^2 \\
\notag & \quad + m(m-1)(2-k)u^2\lv Lv\rv^2 \\
\notag & \quad + 2(k-1)m(m+n-2)uv\lp Lu,Lv\rp \\
\label{eqn:mu} (m+n-2)k\mu u^k v^{2-k} & = -m(m+n-2)uv\lp Lu,Lv\rp + m(m-1)u^2\lv Lv\rv^2 .
\end{align}
\item If also $w$ is a critical point of the conformal GNS functional for variations in $g$, then it holds that
\begin{equation}
\label{eqn:qe}
\begin{split}
0 & = (m+n-2)v\nabla(Lu) - mu\nabla(Lv) \\
\notag & \quad - \frac{1}{n}\left((m+n-2)\lp\nabla(Lu),Lv\rp-m\lp\nabla(Lv),Lu\rp\right)X .
\end{split}
\end{equation}
\end{enumerate}
\end{prop}

\begin{remark}
\label{rk:other_vars}
If one only assumes that $u$ corresponds to a critical point of a conformal GNS functional for variations of $w$, then it holds that
\begin{equation}
\label{eqn:tractor_general_crit_point}
\begin{split}
& \qquad (m+n)\lambda v^2 - 2(m+n-k)\mu u^{k}v^{2-k} \\
& = -(m+n)(m+n-1)v^2\lv Lu\rv^2 + 2m(m+n-1)uv\lp Lu,Lv\rp \\
& \quad - m(m-1)u^2\lv Lv\rv^2 .
\end{split}
\end{equation}
\end{remark}

\begin{remark}
\label{rk:dd_proof}
Note that~\eqref{eqn:tractor_general_crit_point} easily verifies that the functions~\eqref{eqn:gns_dd_minimizer} are critical functions for the corresponding conformal GNS functional.  Indeed, for the standard Euclidean space $(\bR^n,g,1^m\dvol)$, we have that $\lv Lv\rv^2=0$, while the the GNS inequality~\eqref{eqn:gns_dd} corresponds to finding $\sigma_{1,1}$.  Thus, if one can find $u\in\mE[1]$ such that $\lv Lu\rv^2$ and $\lp Lu,Lv\rp$ are negative constants, then~\eqref{eqn:tractor_general_crit_point} holds with $\lambda,\mu>0$.  Using the properties of the basis~\eqref{eqn:rn_tractor_basis}, it is easy to check that this holds if and only if there are constants $a_0,\dotsc,a_n,a_{n+1}$ such that
\[ u = a_0\frac{1+r^2}{2} + \sum_{i=1}^n a_i x^i + a_{n+1}\frac{1-r^2}{2}, \]
$a_0>0$, and $a_0^2-a_{n+1}^2>\sum_{i=1}^n a_i^2$, giving the entire family~\eqref{eqn:gns_dd_minimizer} of minimizers.  To show that these are the only such solutions to~\eqref{eqn:tractor_general_crit_point}, and thus extremal functions, one can proceed as in~\cite{DelPinoDolbeault2002} by invoking results of Pucci and Serrin~\cite{PucciSerrin1998} and of Serrin and Tang~\cite{SerrinTang2000}.
\end{remark}

\begin{proof}

From Proposition~\ref{prop:variational}, it follows easily that
\begin{align*}
nk(m+n-2)\lambda u^{-2} & = k(m+n-2)\left(R_{f,\phi}^m-m\Delta_\rho\beta\right) \\
& \quad + m(2-k)\left(R_{f,\phi}^m-(m+n)\Delta_\rho\beta\right) \\
nk(m+n-2)\mu u^{k-2} v^{-k} & = m\left(R_{f,\phi}^m - (m+n)\Delta_\rho\beta\right) \\
0 & = (\Ric_{f,\phi}^m)_0 .
\end{align*}
The result then follows from Proposition~\ref{prop:smms_tractor_formulae}.
\end{proof}

With Proposition~\ref{prop:tractor} in hand, we are now in a position to classify the GNS constants whose minimizers are also critical for variations in the metric or the measure, under the additional assumption that the underlying smooth metric measure space $(M^n,g,v^m\dvol)$ is such that $v^{-2}g$ is an Einstein metric, or equivalently, that $Lv$ is parallel.  First, if the minimizers are also critical for variations of the metric, then there are three possibilities.

\begin{thm}
\label{thm:qe_rigid}
Let $(M^n,g,v^m\dvol)$ be a smooth metric measure space such that $Lv\in\mT$ is parallel.  Suppose that $w=u^{-\frac{m+n-2}{2}}$ is a critical point of the conformal GNS functional for variations of $g$ and $w$, and moreover, assume that $\mQ_k(w)>0$.  Then $Lu\in\mT$ is parallel.  Moreover, $\lv Lu\rv^2<0$, and one of the three mutually exclusive statements is true.
\begin{enumerate}
\item $\frac{u}{v}$ is constant; in particular, $\lv Lv\rv^2<0$.
\item $k=1$ and $\lv Lv\rv^2=0$.
\item $k=2$ and $\lp Lu,Lv\rp=0$; in particular, $\lv Lv\rv^2>0$.
\end{enumerate}
\end{thm}

\begin{remark}

In tensorial language, Theorem~\ref{thm:qe_rigid} states that if $(M^n,g,1^m\dvol)$ is a smooth metric measure space such that $g$ is Einstein, the GNS constant satisfies $\sigma_{1,k}>0$, and an extremal function $w$ for $\sigma_{1,k}$ is also a critical point of the conformal GNS functional for variations of the metric and the measure, then the metric $\hat g=w^{\frac{4}{m+n-2}}g$ is Einstein with positive scalar curvature and one of the three mutually exclusive statements is true.
\begin{enumerate}
\item $w$ is constant, and in particular, $g$ has positive scalar curvature.
\item $k=1$ and $g$ is Ricci flat.
\item $k=2$ and $g$ has negative scalar curvature.
\end{enumerate}
This is why we say that $\sigma_{1,2}$ is better behaved on hyperbolic space than on Euclidean space.
\end{remark}

\begin{proof}

If $\nabla(Lv)=0$, then it follows from~\eqref{eqn:qe} that $\nabla Lu=\alpha\otimes X$ for some one-form $\alpha$.  By Lemma~\ref{lem:prolongation}, we thus have that $\nabla(Lu)=0$.  Thus the inner products $\lv Lu\rv^2$, $\lp Lu,Lv\rp$, and $\lv Lv\rv^2$ are all constant.  By Lemma~\ref{lem:critical_g} and~\eqref{eqn:lambda}, we then see that
\begin{equation}
\label{eqn:parallel_quadratic}
a\left(\frac{u}{v}\right)^2 + 2b\left(\frac{u}{v}\right) + c = 0 ,
\end{equation}
with the constants $a,b,c$ given by
\begin{equation}
\label{eqn:parallel_quadratic_constants}
\begin{split}
a & = m(m-1)(2-k)\lv Lv\rv^2 \\
b & = m(m+n-2)(k-1)\lp Lu,Lv\rp \\
c & = (m+n-2)k\left(\lambda + (m+n-1)\lv Lu\rv^2\right) .
\end{split}
\end{equation}
Since $\mQ_k(w)>0$, it follows that $\lambda,\mu>0$.  In particular, if $\frac{u}{v}$ is not constant, then it must be the case that $a=b=c=0$, giving the last two cases of the theorem.  If instead $\frac{u}{v}$ is constant, it follows immediately from~\eqref{eqn:lambda} and~\eqref{eqn:mu} that $\lv Lu\rv^2<0$, as desired.
\end{proof}

\begin{remark}

As the proof makes clear, the positivity assumption $\mQ_k(w)>0$ is only needed to establish the signs of the various inner products between $Lu$ and $Lv$.
\end{remark}

Let us now suppose additionally that $(M^n,g,v^m\dvol)$ is such that $v^{-2}g$ is Ricci flat, with the goal of understanding minimizers for GNS constants which are also critical for variations in the measure.  By Theorem~\ref{thm:qe_rigid}, we know that unless $k=1$, the minimizers of the conformal GNS functional will not be critical points for variations in the metric $g$.  We show that they will not even be minimizers for variations in the density $v=1$, thus establishing Theorem~\ref{thm:main_thm} as a corollary.

\begin{thm}
\label{thm:v_rigid}
Let $(M^n,g,v^m\dvol)$ be a smooth metric measure space such that $v^{-2}g$ is Ricci flat and suppose that $w=u^{-\frac{m+n-2}{2}}$ is a critical point of the conformal GNS functional $\mQ_k$ with respect to variations of $w$ and $v$, and moreover suppose that $\mQ_k(w)>0$.  Then $k=1$ and $Lu$ is parallel.
\end{thm}

\begin{proof}

By the assumption on $v$, we may apply Proposition~\ref{prop:tractor} to see that
\begin{equation}
\label{eqn:v_rigid_norms}
\begin{split}
-2(k-1)\mu\left(\frac{u}{v}\right)^k & = \lambda + (m+n-1)\lv Lu\rv^2 \\
-k\mu\left(\frac{u}{v}\right)^{k-1} & = m\lp Lu,Lv\rp .
\end{split}
\end{equation}
In particular, if $k=1$, we see that both $\lv Lu\rv^2$ and $\lp Lu,Lv\rp$ are constant, at which point Lemma~\ref{lem:tractor_lemma} immediately yields that $Lu$ is parallel.  In fact, Lemma~\ref{lem:tractor_lemma} and~\eqref{eqn:v_rigid_norms} together also imply that $k=1$.  To see this, first observe that
\begin{align*}
\nabla\lv Lu\rv^2 & = -\frac{2k(k-1)\mu}{m+n-1}\left(\frac{u}{v}\right)^{k-1}\nabla\frac{u}{v} \\
\nabla\lp Lu,Lv\rp & = -\frac{k(k-1)\mu}{m}\left(\frac{u}{v}\right)^{k-2}\nabla\frac{u}{v} .
\end{align*}
For brevity, denote $C=\frac{k(k-1)\mu}{m(m+n-1)}$ and $X=v\nabla u-u\nabla v$.  Then Lemma~\ref{lem:tractor_lemma} yields
\begin{equation}
\label{eqn:divergence_form_eqn}
\begin{split}
-\frac{v}{u}\lv\nabla Lu\rv^2 & = C\bigg[(m+n-1)u^{n-2}\delta\left(u^{k-n}v^{-k}X\right) \\
& \qquad - mu^{-1}v^{n-1}\delta\left(u^{k-1}v^{1-k-n}X\right)\bigg] \\
& = C\left[\delta X + (k-m-n)u^{-1}\lp X,\nabla u\rp - (k-m)v^{-1}\lp X,\nabla v\rp\right] \\
& = Cu^{m+n-2}v^{-m}\delta\left(u^{k-m-n}v^{m-k}X\right) .
\end{split}
\end{equation}
Thus $\lv\nabla Lu\rv^2$ is a pure divergence.

Finally, the assumption that $u$ is a minimizer of the conformal GNS functional implicitly includes the assumptions that
\[ \int_M \lv\nabla u\rv^2 u^{-m-n}v^m\dvol, \int_M \lv\nabla v\rv^2 u^{2-m-n}v^m\dvol, \int_M u^{2-m-n}v^{m-2}\dvol < \infty . \]
Integrating~\eqref{eqn:divergence_form_eqn} with respect to $u^{2-m-n}v^m\left(\frac{u}{v}\right)^{1-k}\dvol$ then implies that $\lv\nabla Lu\rv^2$ vanishes.  In particular, $Lu$ is parallel, and hence by Theorem~\ref{thm:qe_rigid}, $k=1$, as desired.
\end{proof}
\section{Concluding Remarks}
\label{sec:conclusion}

Let us conclude by making precise how to use our methods to understand GNS inequalities~\eqref{eqn:gns} in the full range $1\leq p\leq q\leq\frac{2n}{n-2}$.  In particular, we will also give a variational characterization of the second family~\eqref{eqn:gns_dd_other} of GNS inequalities whose sharp constants were computed by Del Pino and Dolbeault~\cite{DelPinoDolbeault2002}.

The main observation is that in our definitions of a smooth metric measure space and geometric notions associated to them, the constraint $m\geq 0$ is not actually needed; we only require that $m+n-2\not=0$.  An important special case is when we take $m\in[-\infty,2-n)$.  For example, it is this perspective which Chen~\cite{Chen2007} used to construct new examples of conformally compact Einstein manifolds, and which the author~\cite{Case2010b} used to prove a precompactness theorem for quasi-Einstein manifolds.

Recall our definition of the GNS constants
\[ \sigma_{1,k}\left(g,e^{-\phi}\dvol,m\right) = \inf \left\{ \frac{\mF(w) \left(\int_M w^{\frac{2(m+n-k)}{m+n-2}}v^{m-k}\right)^p}{\left(\int w^{\frac{2(m+n)}{m+n-2}}v^m\dvol\right)^q} \right\}, \]
where $p=\frac{2m}{nk}$, $q=\frac{2m+k(n-2)}{nk}$, and the infimum is taken over all nonzero $w\in L_1^2(M)$.  This makes sense if we instead suppose that $m\leq-n-2$ and $k\in[0,-\frac{2m}{n-2}]$.  Indeed, on $\bR^n$ for this range of $m$ and $k$, $\sigma_{1,k}$ is equivalent to the sharp constant in the GNS inequality
\begin{equation}
\label{eqn:gns_constant_other}
\lV w\rV_{\frac{2(m+n-k)}{m+n-2}} \leq C \lV\nabla w\rV_2^\theta \lV w\rV_{\frac{2(m+n)}{m+n-2}}^{1-\theta} .
\end{equation}
In other words, combined with Lemma~\ref{lem:gns_const} for the case $m\geq 0$, we see that the conformal GNS constants do include all cases of the sharp constants in the GNS inequalities~\eqref{eqn:gns}.  Moreover, just like in the case $m\in[0,\infty]$, taking $k=1$ in~\eqref{eqn:gns_constant_other} recovers the second family~\eqref{eqn:gns_dd_other} studied by Del Pino and Dolbeault~\cite{DelPinoDolbeault2002}.

Now, a quick check of the proofs of Proposition~\ref{prop:smms_tractor_formulae}, Theorem~\ref{thm:qe_rigid}, and Theorem~\ref{thm:v_rigid} reveals that the assumption $m\geq 0$ is unnecessary; one need only be careful in considering the relationship between the signs of $\lambda$ and $\mu$ and of the inner products involving $Lu$ and $Lv$.  In particular, a straightforward modification of Remark~\ref{rk:dd_proof} shows that in the case $k=1$, extremal functions for the family~\eqref{eqn:gns_dd_other} of GNS inequalities can be constructed from the function
\[ w(x) = \left(1-\lv x\rv^2\right)_+^{-\frac{m+n-2}{2}}, \]
which is the conformal factor of the hyperbolic metric on the ball $B(1)$, and we again have the characterization of the family~\eqref{eqn:gns_dd_other} as the only family of GNS inequalities of the form~\eqref{eqn:gns_constant_other} whose extremal functions are also critical points for variations in the metric or the measure.

\bibliographystyle{abbrv}
\bibliography{../bib}
\end{document}